\documentclass[11pt]{amsart}
\usepackage{xypic}

\newcommand{\Bred}{\lfloor B \rfloor}

\newcommand{\llc}{\lfloor}
\newcommand{\lrc}{\rfloor}

\newcommand{\mb}{\mathbb}
\newcommand{\mc}{\mathcal}
\newcommand{\wt}{\widetilde}
\newcommand{\wX}{\widetilde{X}}
\newcommand{\wB}{\widetilde{B}}
\newcommand{\uo}{^{(1)}}
\newcommand{\ut}{^{(2)}}

\newcommand{\inv}{^{-1}}
\newcommand{\sat}{^{\rm sat}}
\newcommand{\sn}{^{\rm sn}}

\newcommand{\red}{_{\rm red}}
\newcommand\cO{{\mathcal O}}
\newcommand\bP{{\mathbb P}}
\newcommand\bQ{{\mathbb Q}}

\newcommand\bF{{\mathbb F}}

\DeclareMathOperator{\Exc}{\operatorname{Exc}}
\DeclareMathOperator{\LCS}{\operatorname{LCS}}

\DeclareMathOperator{\Supp}{\operatorname{Supp}}
\DeclareMathOperator{\Spec}{\operatorname{Spec}}

\DeclareMathOperator{\codim}{\operatorname{codim}}

\DeclareMathSymbol{\twoheadrightarrow}  {\mathrel}{AMSa}{"10}

\newtheorem{theorem}{Theorem}[section]
\newtheorem{lemma}[theorem]{Lemma}

\newtheorem{corollary}[theorem]{Corollary}
\theoremstyle{definition}
\newtheorem{definition}[theorem]{Definition}
\newtheorem{discussion}[theorem]{Discussion}

\newtheorem{example}[theorem]{Example}

\newtheorem{remark}[theorem]{Remark}
\newtheorem*{acknowledgements}{Acknowledgements}

\newcommand\defn[1]{{\bf #1}}

\title{Limits of stable pairs}
\author{Valery Alexeev}
\address{Department of Mathematics, University of Georgia, 
  Athens GA 30602, USA}
\email{valery@math.uga.edu}
\date{July 21, 2006; revised September 29, 2007}

\begin{document}
\maketitle
\begin{abstract}
  Let $(X_0,B_0)$ be the canonical limit of a one-parameter family of stable
  pairs, provided by the log Minimal Model Program. We prove that
  $X_0$ is $S_2$ and that $\lfloor B_0 \rfloor$ is $S_1$, as an application
  of a general local statement: if $(X,B+\epsilon D)$ is
  log canonical and $D$ is $\bQ$-Cartier then $D$ is $S_2$ and
  $\lfloor B \rfloor \cap D$ is $S_1$, i.e. has no embedded components.

  When $B$ has coefficients $<1$, examples due to Hacking and Hassett
  show that $B_0$ may indeed have embedded primes. We resolve this
  problem by introducing a category of stable branchpairs. We prove
  that the corresponding moduli functor is proper for families with
  normal generic fiber.
\end{abstract}
\tableofcontents

Let $U=S\setminus 0$ be a punctured nonsingular curve, and
$f:(X_U,B_U)\to Y\times U$ be a family of stable maps (precise
definitions follow).  It is well understood, see e.g.
\cite{KollarShepherdBarron,Alexeev_LCanModuli} that log Minimal Model
Program leads to a natural completion of this family over $S$,
possibly after a finite ramified base change $S'\to S$.  This, in
turn, leads to the construction of a proper moduli space of stable
maps (called stable pairs if $Y$ is a point) once some standard
conjectures, such as log MMP in dimension $\dim X+1$ and boundedness,
and some technical questions have been resolved.

The purpose of this paper is solve two such technical issues. The
first one is the Serre's $S_2$-property for the one-parameter limits,
which implies that the limit is semi log canonical:

\begin{theorem}\label{thm:main}
  Let $(X,B)\to S$ be the stable log canonical completion of a family of
  log canonical pairs. Then for the central fiber one has:
  \begin{enumerate}
  \item $X_0$ is $S_2$,
  \item $\llc B_0 \lrc$ is $S_1$, i.e. this scheme has no embedded components. 
  \end{enumerate}
\end{theorem}
As a corollary, if $B$ is reduced (i.e. all $b_j=1$) then the central
fiber $f:(X_0,B_0)\to Y$ is a stable map. 

For surfaces with reduced $B$, Theorem~\ref{thm:main} was proved by
Hassett \cite{Hassett_StableLimits}. Even in that case, our proof is
different. Whereas the proof in \cite{Hassett_StableLimits} is global,
i.e. it requires an actual semistable family of projective surfaces
with relatively ample $K_X+B$, our proof is based on the following
quite general local statement:

\begin{lemma}
  Let $(X,B)$ be a log canonical pair which has no zerodimensional
  centers of log canonical singularities. Then for every closed point
  $x\in X$, the local ring $\cO_{X,x}$ is $S_3$.
\end{lemma}


As a consequence, we obtain the following theorem from which
\eqref{thm:main} follows at once. 

\begin{theorem}
  Let $(X,B)$ be a log canonical pair and $D$ be an effective Cartier
  divisor. Assume that for some $\epsilon>0$ the pair $(X,B+\epsilon
  D)$ is log canonical. Then $D$ is $S_2$ and $\Bred\cap D$ is $S_1$. 
\end{theorem}

The second question we consider is the following. When the coefficients
$b_j$ in are less than one, Hacking and Hassett gave
examples of families of stable surface pairs in which
the central fiber $B_0$ of $B$ 
 indeed does have embedded primes. We
resolve this problem by introducing, following ideas of
\cite{AlexeevKnutson}, a new category, that of \emph{stable 
  branchpairs}, which avoids nonreduced schemes. We define the
moduli functor in this category and check the valuative criterion of
properness for families with normal generic fiber.


With branchdivisors thus well-motivated, we define, in a
straighforward way,  branchcycles of other dimensions as well. 

\begin{acknowledgements}
  It is a pleasure to acknowledge helpful conversations with Florin
  Ambro, Paul Hacking, Brendan Hassett, J\'anos Koll\'ar
  and Allen Knutson. I also thank the referee for thoughtful comments.
  The research was partially supported by NSF under the grant DMS-0401795.
\end{acknowledgements}

Throughout most of the paper we work over an algebraically closed field of
characteristic zero, and relax this condition to an arbitrary field
for the last section.

\section{Basic definitions}
\label{sec:basic-definitions}

All \defn{varieties} in this paper will be assumed to be connected and
reduced but not necessarily irreducible.  A \defn{polarized variety}
is a projective variety $X$ with an ample invertible sheaf~$L$.
A pair $(X,B)$ will always consist of a variety $X$ and a
$\bQ$-divisor $B=\sum b_j B_j$, where $B_j$ are effective Weil
divisors on $X$, and $0< b_j\le 1$. 

We use standard definitions and notations of Minimal Model Program
for discrepancies $a(X,B,E_i)$, notions of log canonical pairs
(abbreviated lc), klt pairs, etc., as in \cite{KollarMori_Book}.
We assume standard definitions from commutative algebra for
the Serre's conditions $S_n$. We now list the slightly less standard
definitions. 

\begin{definition}
  Let $(X,B)$ be an lc pair. A \defn{center of log canonical
    singularities} of $(X,B)$ (abbreviated to a center of $\LCS(X,B)$)
  is the image of a divisor $E_i\subset Y$ on a resolution $f:Y\to X$
  that has discrepancy $a(X,B,E_i)=~-1$. 

  If $f:Y\to X$ is log a smooth resolution of $(X,B)$ and $E=\sum E_i$
  is the union of all divisors with discrepancy $-1$ (some
  exceptional, some strict preimages of components of $B$ with
  $b_j=1$) then the centers of $\LCS(X,B)$ are the images of the
  nonempty strata $\cap E_i$.
\end{definition}

We will use the following important results of Florin Ambro,
which were further clarified by Osamu Fujino. The first is Ambro's
generalization of Koll\'ar's injectivity theorem
\cite{Kollar_DirectImages1}, and the second describes properties of
log centers. 

\begin{theorem}[Injectivity for varieties with normal crossings, simple form]
\label{thm:Ambro_vanishing}
  Let $Y$ be a nonsingular variety, $E+S+\Delta$ be a normal crossing
  $\mb R$-divisor on $Y$, $E,S$ and $\Delta$ have no components in
  common,  $E+S$ is reduced, and $\llc\Delta\lrc=0$. 
  
  Let $f:Y\to X$ be a proper morphism, $A$ a Cartier divisor on $E$,
  and assume that the divisor $H\sim_{\mathbb R} A - (K_E+S+\Delta)$ on
  $E$ is $f$-semiample.  Then every nonzero section of $R^i f_* \mathcal
  O_E(A)$ contains in its support the $f$-image of some strata of
  $(E,S+\Delta)$. 

  Here, $K_E$ stands for the dualizing invertible sheaf $\omega_E$,
  and the strata of $(E,S+\Delta)$ are the intersections of the
  components of $E$ and $S$.
\end{theorem}
\begin{proof}
  This is a special case of \cite[3.2(i)]{Ambro_QLvars}, see also
  \cite{Ambro_LogCenters} for another exposition. This theorem was
  also reproved  in \cite[5.7,5.15]{Fujino_VanishingInjectivity}, see
  also \cite{Fujino_NotesLMMP}.
\end{proof}

\begin{theorem}[Properties of log centers]\label{thm:centers}
  \begin{enumerate}
  \item Every irreducible component of the intersection of two centers
    is a center.
  \item For any $x\in X$ the minimal center containing $x$ is normal. 
  \item A union of any set of centers is seminormal.
  \end{enumerate}
\end{theorem}
\begin{proof}
  (1) and (2) are contained in \cite{Kawamata_FujitasFreeness,
    Kawamata_Subadjunction2} in the case when there exists a klt
  pair $(X,B')$ with $B'\le B$. For  the general case these are in 
  \cite[4.8]{Ambro_QLvars}.
  (3) is \cite{Ambro_LCS} and \cite[4.2(ii),4.4(i)]{Ambro_QLvars}.

  Also, a very easy, one-page proof of these properties, which uses
  only the above injectivity theorem, is contained in
  \cite[\S4]{Ambro_LogCenters}. 
\end{proof}

\begin{definition}
  A pair $(X,B)$ is called \defn{semi log canonical} (slc) if 
  \begin{enumerate}
  \item $X$ satisfies Serre's condition $S_2$,
  \item $X$ has at worst double normal crossing singularities in
    codimension one, and 
    no divisor $B_j$ contains any component of this double locus,
  \item some multiple of the Weil $\mathbb Q$-divisor $K_X+B$, well
    defined thanks to the previous condition, is $\mathbb Q$-Cartier, and
  \item denoting by $\nu:X^{\nu}\to X$ the normalization, the pair
    \linebreak
    $(X^{\nu},\ \text{(double locus)} + \nu\inv B )$ is log canonical.
  \end{enumerate}
\end{definition}

\begin{definition}\label{defn:stable-pair}
A pair  $(X,B=\sum b_j B_j)$ (resp. a map $f:(X,B)\to Y$) is called
  a \defn{stable map} if the following two conditions are
  satisfied: 
  \begin{enumerate}
  \item \emph{on singularities:} the pair $(X,B)$ is semi log canonical, and
  \item \emph{numerical:} the divisor $K_X+B$ is ample (resp. $f$-ample).
  \end{enumerate}
  A \defn{stable pair} is a stable map to a point. 
\end{definition}

\begin{definition}
  A variety $X$ is \defn{seminormal} if any proper bijection $X'\to X$
  is an isomorphism.
\end{definition}
It is well-known, see e.g. \cite[I.7]{Kollar_RationalCurves},
that every variety has a unique seminormalization
$X\sn$ and it has a universal property: any morphism $Y\to X$ with
seminormal $Y$ factors through $X\sn$.

\section{$S_2$ and seminormality}
\label{sec:S2-seminormality}

We collect some mostly well-known facts about the way
the $S_2$ property and seminormality are related.

\begin{definition}
  The \defn{$S_2$-fication}, or \defn{saturation in codimension 2}
  of a variety $X$ is defined to be 
  \begin{displaymath}
    \pi\sat_X: X\sat = \varinjlim\  \Spec_{\cO_{X}} \cO_{X\setminus Z}  
    \ \to \ X
  \end{displaymath}
  in which the limit goes over closed subsets  $Z\subset X$ with
  $\codim_X Z\ge2$.
  The morphism $\pi\sat_X$ is finite: indeed, it is dominated by the
  normalization of $Y$. 

  More generally, for any closed subset $D\subset X$ the
  \defn{saturation in codimension 2 along $D$}
  \begin{displaymath}
    \pi\sat_{X,D}:X\sat_D  \to X
  \end{displaymath}
  is defined by taking the limit as above that goes only over
  $Z\subset  D$. Hence, $\pi\sat_X = \pi\sat_{X,X}$. 
\end{definition}

\begin{lemma}\label{lem:sat-along-D-S2}
  $\pi\sat_{X,D}$ is an isomorphism iff for any subvariety $Z\subset
  D$ the local ring $\cO_{X,Z}$ is $S_2$. 
\end{lemma}
\begin{proof}
  Let $Z\subset D$ be a subvariety with $\codim_X Z\ge 2$. 
  By the cohomological characterization of
  depth (see f.e.  \cite[Thm.  28]{Matsumura_ComRingThry} or
  \cite[18.4]{Eisenbud_CommAlg}) the local ring $\mc O_{X,Z}$ has
  depth $\ge2$ iff any short exact sequence 
  \begin{displaymath}
    0\to \mc O_{X,Z} \to F \to Q \to 0
  \end{displaymath}
  of $\cO_{X,Z}$-modules with $\Supp Q=Z$ splits. 

  If $\pi\sat_D$ is an isomorphism then for every exact sequence as
  above $F\sat_D =\mc O_{X,Z}$, and the canonical restriction morphism
  $F\to F\sat_D$ provides the splitting.
  If $\pi\sat_D$ is not an isomorphism over some $Z\subset D$ then the
  localization of
  \begin{displaymath}
    0\to O_X \to \pi\sat_* \cO_{X\sat} \to Q \to 0
  \end{displaymath}
  at $Z$ does not split and $Q\ne0$.
\end{proof}

\begin{lemma}\label{lem:seminormal=>S2}
  Assume that $X$ is seminormal and $\pi\sat_{X,D}$ is a bijection.
  Then for any subvariety $Z\subset D$ the local ring $\cO_{X,Z}$ is
  $S_2$.
\end{lemma}
\begin{proof}
  Since $X$ is seminormal, $\pi\sat_{X,D}$ is an isomorphism, so the
  the previous lemma applies.
\end{proof}

\begin{lemma}\label{lem:S2=>seminormal}
  Assume $X$ is $S_2$ and is seminormal in codimension 1. Then $X$ is
  seminormal. 
\end{lemma}
\begin{proof}
  We have $(X\sn)\sat = X\sat$ and $X\sat=X$, hence $X\sn \to
  X$ is an isomorphism.
\end{proof}

\begin{corollary}
  Semi log canonical $\implies$ seminormal.
\end{corollary}

\section{Singularity theorems}
\label{sec:singularity-theorems}

Let $X$ be a normal variety, which by Serre's criterion implies that
$X$ is $S_2$. Let $f:Y\to X$ be a resolution of singularities. Then we
have:

\begin{lemma}
  Assume $\dim X>2$. Then $X$ is $S_3$ at every closed point $x\in X$ iff
  $R^1f_*\cO_Y$ has no associated components of dimension 0, i.e. the
  support of every section of $R^1f_*\cO_Y$ has dimension $>0$.
\end{lemma}
\begin{proof}
  By considering an open affine neighborhood of $x$ and then
  compactifying, we can assume that $X$ is projective with an ample
  invertible sheaf $L$. (Since the property of being $S_3$ at closed
  points is open, one can compactify without introducing ``worse''
  points.)  Then by the proof of
  \cite[Thm.III.7.6]{Hartshorne_AlgebraicGeometry}, $X$ is $S_3$ at every
  closed point iff for all $r\gg0$ one has $H^2\big( \cO_X(-rL)
  \big)=0$.

  The spectral sequence
  \begin{displaymath}
    E_2^{p,q} = H^p\big( R^qf_* \cO_Y(-rL) \big)
    \Rightarrow H^{p+q}\big( \cO_Y (-rf^*L) \big)
  \end{displaymath}
  together with the fact that $H^{1\text{ and } 2}\big( \cO_Y (-rf^*L)
  \big)=0$ by Generalized Kodaira's vanishing theorem
  \cite[2.70]{KollarMori_Book}, imply that
  \begin{displaymath}
    d_2^{0,1}: H^0\big( R^1f_*\cO_Y (-rL) \big) 
    \to H^2\big( \cO_X(-rL) \big)
  \end{displaymath}
  is an isomorphism. Further, $H^0\big( R^1f_*\cO_Y (-rL)
  \big)=0$ for $r\gg0$ precisely when the sheaf $R^1f_*\cO_Y$ has no
  associated components of dimension~0.
\end{proof}

Log terminal pairs have rational singularities, and hence are
Cohen-Maca\-ulay, see \cite[Thm.5.22]{KollarMori_Book} for a simple
proof. Log canonical singularities need not be $S_3$. The easiest
example is a cone over an abelian surface $S$. Indeed, in this case
$R^1f_*\cO_Y = H^1\big( \cO_S \big)$ is non-zero and
supported at one point. However, we will prove the following:

\begin{lemma}\label{lem:S3}
  Let $(X,B)$ be a log canonical pair which has no zerodimensional
  centers of log canonical singularities. Then for every closed point
  $x\in X$ the local ring $\cO_{X,x}$ is $S_3$.
\end{lemma}
\begin{proof}
As in the previous proof, we can assume that $(X,L)$ is a polarized variety,
and we must prove that for $r\gg0$
one has $H^2\big(\cO_X(-rL)\big)=0$. 
Let $f:Y\to X$ be a resolution of singularities of $(X,B)$ such that
$f\inv B\cup \Exc(f)$ is a divisor with global normal crossings. Then we
can write
\begin{displaymath}
  K_Y \sim_{\bQ} f^*(K_X+B) - E + A - \Delta,
\end{displaymath}
where
\begin{enumerate}
\item $E=\sum E_j$ is the sum of the divisors $B_j$ with $b_j=1$ and
  the exceptional divisors of $f$ with discrepancy $-1$,
\item $A$ is effective and integral,
\item $\Delta$ is effective and $\llc\Delta\lrc =0$.
\end{enumerate}

Since the pair $(X,B)$ is log canonical and the coefficients of $B$
satisfy $0<b_j\le 1$, it follows that $A$ is $f$-exceptional, 
$E$ has no components in common with $\Supp\Delta$ and with $A$,
and the union $E\cup\Supp A\cup\Supp\Delta$ is a divisor with global
normal crossings.

Then $-E+A \sim_{\bQ} K_Y + \Delta - f^*(K_X+B)$. The Generalized
Kodaira Theorem (Kawamata-Viehweg theorem) gives $R^qf_*
\cO_Y(-E+A)~=~0$ for $q>0$. Therefore, by pushing forward the exact
sequence
\begin{displaymath}
  0 \to \cO_Y(-E+A) \to \cO_Y(A) \to \cO_E(A) \to 0
\end{displaymath}
we obtain $R^1f_* \cO_Y(A) \simeq R^1f_* \cO_E(A)$. Now, on $E$ one has
\begin{displaymath}
  A \sim_{\bQ} K_E + \Delta - f^*(K_X+B),
\end{displaymath}
where $K_E$ stands for the (invertible) dualizing sheaf $\omega_E$.
Therefore, by Ambro's injectivity theorem \ref{thm:Ambro_vanishing},
applied here with $H=-f^*(K_X+B)$ and $S=0$,
 the
support of every nonzero 
section of the sheaf $R^1f_*\cO_E(A)$ contains a center of
$\LCS(X,B)$, hence has dimension $>0$. 

Now consider the following commutative diagram
\begin{displaymath}
  \diagram
  H^2\big( \cO_Y(-rf^*L) \big) \rto
&  H^2\big( \cO_Y(A-rf^*L) \big) \\
  H^2\big( f_*\cO_Y(-rf^*L) \big) \uto \ar@{=}[d]
&  H^2\big( f_*\cO_Y(A-rf^*L) \big) \uto_{\phi} \ar@{=}[d] \\
  H^2\big( \cO_X(-rL) \big)  \ar@{=}[r]
&  H^2\big( \cO_X(-rL) \big) 
  \enddiagram
\end{displaymath}
Since by Generalized Kodaira's vanishing theorem $H^2\big(
\cO_Y(-rf^*L) \big)=0$, this implies $H^2\big( \cO_X(-rL) \big) =0$ if
we could prove that $\phi$ is injective. Finally, the spectral
sequence
  \begin{displaymath}
    E_2^{p,q} = H^p\big( R^qf_* \cO_Y(A)(-rL) \big)
    \Rightarrow E^{p+q} = H^{p+q}\big( \cO_Y (A-rf^*L) \big)
  \end{displaymath}
  in a standard way produces the exact sequence
  \begin{displaymath}
    E_2^{0,1} \overset{d_2^{0,1}}{\longrightarrow}
    E_2^{2,0} \longrightarrow E^2
  \end{displaymath}
  In our case, $E_2^{0,1}=H^0\big ( R^1f_* \cO_Y(A)(-rL) \big)=0$ by
  what we proved above (the sheaf $R^1f_* \cO_Y(A)$ has no associated
  components of dimension 0), and the second homomorphism is
  $\phi$. Hence, $\phi$ is injective. This completes the proof. 
\end{proof}

\begin{remark}
  One has to be careful that \eqref{lem:S3} does \emph{not} imply
  that $X$ is $S_3$. Indeed, let $X'$ be a variety which is $S_2$ but not
  $S_3$, for example a cone over an abelian surface, and let $X$ be the
  cartesian product of $X'$ with a curve $C$. Then $X$ is $S_3$ at every
  closed point but not at the scheme point corresponding to
  (vertex)$\times C$. 
\end{remark}

\begin{theorem}\label{thm:S2}
  Let $(X,B)$ be an lc pair and $D$ be an effective Cartier
  divisor. Assume that for some $\epsilon>0$ the pair $(X,B+\epsilon
  D)$ is lc. Then $D$ is $S_2$.
\end{theorem}
\begin{proof}
  Suppose that for some subvariety $Z\subset D$ the local ring
  $\cO_{D,Z}$ is not $S_2$, then $\cO_{X,Z}$ is not $S_3$. Let
  \begin{displaymath}
    (X^{(d)},B^{(d)}) = (X,B) \cap H_1 \cap \dots \cap H_d
  \end{displaymath}
  be the intersection with $d=\dim Z$ general
  hyperplanes such that $Z^{(d)}=Z\cap H_1\cap\dots\cap H_d\ne
  \emptyset$. Then 
  \begin{enumerate}
  \item the pair $(X^{(d)},B^{(d)})$ is lc by the general properties
    of lc (apply Bertini theorem to a resolution), and
  \item   $X^{(d)}$ is not $S_3$ at a closed point $P\in Z^{(d)}$ (by
    the semicontinuity of depth along $Z$ on fibers a 
    morphism; applied to a generic projection $X\to \bP^d$).
  \end{enumerate}

  Let $W$ be a center of $\LCS(X,B)$. Since $(X,B+\epsilon D)$ is lc,
  $W$ is not contained in $D$. Then the corresponding centers,
  irreducible components of $W^{(d)}= W\cap H_1\cap\dots\cap H_d$ are
  not contained in $D^{(d)}$. Hence, by shrinking a neighborhood of
  $D^{(d)}$ in $X^{(d)}$ we can assume that $(X^{(d)},B^{(d)})$ has no
  zerodimensional centers of $\LCS$. But then $X^{(d)}$ is $S_3$ at $P$
  by \eqref{lem:S3}, a contradiction. 
\end{proof}

\begin{theorem}\label{thm:S1}
  Under the assumptions of \eqref{thm:S2}, the scheme $\Bred\cap D$ is $S_1$. 
\end{theorem}
\begin{proof}
  Note that $\Bred$ is a union of several centers of
  $\LCS(X,B)$. As such, it is seminormal by Theorem~\ref{thm:centers}.

  We claim that the saturation $\pi=\pi\sat_{\Bred,\Bred\cap D}$ of
  $\Bred$ in codimension 2 along $\Bred\cap D$ is a bijection.
  Otherwise, there exists a subvariety $Z\subset \Bred$ intersecting
  $D$ such that $\pi: \pi\inv(Z)\to Z$ is several-to-one along $Z$.
  Then cutting by generic hyperplanes, as above, we obtain a pair
  $(X^{(d)},B^{(d)})$ such that $Z^{(d)}$ is a point $P$ and
  $\Bred^{(d)}$ has several analytic branches intersecting at $P$.

  After going to an \'etale cover, which
  does not change the lc condition, we can assume that $P$ is a
  component of the intersection of two irreducible component of the
  locus of $\LCS(X^{(d)},B^{(d)})$. 

  But then $P$ is a center of $\LCS$ itself, by
  Theorem~\ref{thm:centers}(1). 
  This is not possible, again because
  $(X^{(d)},B^{(d)}+\epsilon D)$ is lc; contradiction.  

  The saturation morphism $\pi:\Bred\sat_{\Bred,\Bred\cap D}
  \to \Bred$ is a bijection, and $\Bred$ is seminormal.
  By Lemma~\ref{lem:seminormal=>S2} this implies that
  $\Bred$ is $S_2$ along any subvariety $Z\subset D$. Therefore
  $\Bred\cap D$ is $S_1$.
\end{proof}

\section{One-parameter limits of stable pairs}
\label{sec:one-parameter-limits}

Let $U=(S,0)$ be a punctured nonsingular curve and let $f_U:(X_U,B_U)\to
Y\times U$ be a family of stable maps with normal $X_U$, so that
$(X_U,B_U)$ is lc.

The stable limit of this family is constructed as follows.  Pick some
extension family $f:(X,B)\to Y\times S$. Take a resolution of
singularities, which introduces some exceptional divisors $E_i$.
Apply the Semistable Reduction Theorem to this resolution
together with the divisors, as in \cite[Thm.7.17]{KollarMori_Book}.
The result is that after a ramified base change $(S',0)\to (S,0)$ we
now have an extended family $\tilde f':(\wt X',\wt B')$ such that $\wt
X'$ is smooth, the central fiber $\wt X'_0$ is a reduced normal
crossing divisor, and, moreover, $\wt X'_0 \cup \Supp\wt B'\cup \wt
E'_i$ is a normal crossing divisor. Let us drop the primes in this
notation for simplicity, and write $X,S,$ etc. instead of $X',S'$ etc.
  
It follows that the pair $(\wt X, \wt B+\wt X_0 + \sum \wt E_i)$
has log canonical singularities and is relatively of general type over
$Y\times S$.  Now let $f:(X,B+X_0)\to Y\times S$ be its log
canonical model, guaranteed by the log Minimal Model Program. The
divisor $K_{X}+B+X_0$ is $f$-ample and the pair $(X,B+X_0)$ has
canonical singularities. 

\begin{theorem}
  The central fiber $X_0$ is $S_2$, and the scheme
  $\Bred\cap X_0$ is $S_1$.
\end{theorem}
\begin{proof}
  Immediate from  \eqref{thm:S2} and \eqref{thm:S1} by taking $D=X_0$ and
  $\epsilon=1$. 
\end{proof}

\section{Branchpairs}
\label{sec:branchpairs}

In \cite{Hassett_WeightedStableCurves} Hassett constructed moduli
spaces of weighted stable curves, i.e. one-dimensional pairs $(X,\sum
b_iB_i)$ with $0<b_i\le1$. It is natural to try to extend this
construction to higher dimensions.

However, in the case of surfaces Hacking and Hassett gave examples of
one-parameter families of pairs $(X,bB)\to S$ with irreducible $B$
such that $B_0$ has an embedded point. Such examples are
constructed by looking at families $(X,B)\to S$ in which $B$ is
\emph{not} $\bQ$-Cartier.  Recall that by the definition of a log
canonical pair $K_X+B$ must be $\bQ$-Cartier but neither $K_X$ nor $B$
are required to be such. The following explicit example was
communicated to me by Brendan Hassett, included here with his gracious
permission.  

\begin{example}\label{exmp:Hassett}
  Let $\bF_n$ denote the Hirzebruch ruled surface with exceptional
  section $s_n$ ($s_n^2=-n$) and fiber $f_n$; in particular
  $\bF_0=\bP^1\times \bP^1$ with two rulings denoted by $f_0$ and $s_0$. 

  Let $l\sim s_0 +2f_0$ be a smooth curve in $\bF_0$ and let $\wt X$
  be the blowup of $\bF_0 \times S$ along $l\times 0$ in the central
  fiber. Then $\wt X_0$ is the union of two irreducible components
  $\wX_0\uo=\bF_0$ and $\wX_0\ut=\bF_4$ intersecting along $l$, and
  $l\sim s_4$ in $\bF_4$.

  Let $\wt B_0= \wB_0\uo \cup \wB_0\ut$ be a curve in the
  central fiber such that $\wt B_0\uo \sim 2s_0$ is the union of
  two generic lines and $\wt B_0\ut\sim 4(s_4+4f_4) + 4f_4$,
  intersecting at 4 points $P_1,P_2,P_3,P_4$. Then $\wt B_0$ is a
  nodal curve of genus 35. Let $\wt B$ be a family of curves obtained
  by smoothing $\wt B_0$. 

  Denote by $f:\wt X\to X$ the morphism blowing down the divisor
  $\wX_0\uo$, and $B=f(\wB)$. One easily computes that $K_X$ is not
  $\bQ$-Cartier (because $s_0+2f_0$ is not proportional to
  $K_{\wX_0\uo}$) but $K_X+ 1/2 B$ is, and that $(X,1/2 B)$ has
  canonical singularities. 

  After the blowdown, the curve $(B_0)\red$ in the central fiber is
  obtained from $\wB_0\ut$ by gluing the four points $P_i$ together.
  The curve $\wB_0$ and its smoothings have arithmetic genus 35. The
  curve $(B_0)\red$ has genus 36. Hence, $B_0$ has an embedded point.
\end{example}

So, if one wants to work with arbitrary coefficients, which is very
natural, one must enlarge the category of pairs in some way, or use
some other trick to solve the problem. There are at least two
ways to proceed:

\medskip
(1) One can work with \emph{floating coefficients.} This means that
  we must require the divisors $B_j$ to be $\mb Q$-Cartier, and the
  pairs $(X,\sum (b_j+\epsilon_j)B_j)$ to be semi log canonical and
  ample for all $0<\epsilon_j\ll 1$. Hacking did just that in
  \cite{Hacking_PlaneCurves} for planar pairs $(\mb
  P^2,(3/d+\epsilon)D)$. And the moduli of stable toric,
  resp. abelian pairs in \cite{Alexeev_CMAV} can be interpreted as
  moduli of semi log canonical stable pairs $(X,\Delta+\epsilon B)$,
  resp. $(X,\epsilon B)$.

  However, it is very desirable to work with constant coefficients,
  and the coefficients appearing in the above-mentioned examples are
  fairly simple, such as $b_1=1/2$.

\medskip
(2) One can work with the pairs $(X,\sum b_jB_j)$, where $B_j$ are
  codimension-one subschemes of $X$, possibly with embedded
  components. This can be done in two ways:

(a) \emph{Natural.} One should define (semi) log canonical pairs
    $(X,Y)$ of a variety $X$ with a subscheme $Y$. This was done for
    pairs with \emph{smooth} variety $X$
    (see, e.g. \cite{Mustata_SingsOfPairs}) and more
    generally when $X$ is $\mb Q$-Gorenstein. But: this is
    insufficiently general for our purposes, especially if we
    consider the case of pairs of dimension $\ge 3$.

(b) \emph{Unnatural.} One can work with subschemes $B_j$ that
    possibly have embedded components but then ignore them, by
    saturating in codimension~2. For example, one should define the
    sheaf $\mc O_X( N(K_X+B) )$ as
    \begin{displaymath}
      \mc O_X( N(K_X+B) ) = \varinjlim_{U}\ j_{U*} \mc O_U ( NK_U + B),
    \end{displaymath}
    where the limit goes over open dense subsets $j_U:U\to X$ 
    with $\codim (X\setminus U)\ge2$ such
    that $B\cap U$ has no embedded components and such that $U$ is
    Gorenstein. But this does feel quite artificial.

\bigskip
Building on \cite{AlexeevKnutson}, I now propose a different solution
which avoids nonreduced schemes altogether. 

\begin{definition}
  Let $X$ be a variety of pure dimension $d$. A \defn{prime branchdivisor}
  of $X$ is a variety $B_j$ of pure dimension $d-1$ together with a
  \emph{finite} (so, in particular proper) morphism $\varphi_j:B_j\to X$. 
\end{definition}

Let us emphasize again that by our definition of variety, $B_j$ is
connected, possibly reducible and, most importantly, \emph{reduced}.
Hence, a prime branchdivisor is simply a connected branchvariety, as
defined in \cite{AlexeevKnutson}, of pure codimension $1$.

\begin{definition}
  A \defn{branchdivisor} is an element of a free abelian group
  $bZ_{d-1}(X)$ with prime branchdivisors $B_j$ as generators. If $A$
  is an abelian group (such as $\mb Q$, $\mb R$, etc.) then an
  $A$-branchdivisor is an element of the group $bZ_{d-1}(X)\otimes A$.

  The \defn{shadow} of a branchdivisor $\sum b_jB_j$ is the ordinary
  divisor $\sum b_j \varphi_{j*}(B_j)$ on $X$. We will use the shortcut
  $\varphi_* B$ for the shadow of $B$.
\end{definition}

We will be concerned with $\mb Q$-branchdivisors in this paper,
although $\mb R$-coefficients are frequently useful in other contexts.

\begin{definition}
  A \defn{branchpair} is a pair $(X,\sum b_jB_j)$ of a variety and a $\mb
  Q$-branchdivisor on it, where $B_j$ are prime branchdivisors and
  $0<b_j\le1$. This pair is called \defn{(semi) log
    canonical} (resp. terminal, log terminal, klt) if so is its shadow
  $(X,\sum b_j \varphi_{j*}(B_j))$. 
\end{definition}

\begin{definition}
  A \defn{family of branchpairs} over a scheme $S$ is a 
  morphism $\pi:X\to S$ and finite morphisms $\varphi_j:B_j\to X$ such that
  \begin{enumerate}
  \item $\pi:X\to S$ and all $\pi\circ \varphi_j:B_j\to S$ are flat, and
  \item every geometric fiber $(X_{\bar s}, \sum b_j (B_j)_{\bar s})$
    is a branchpair.
  \end{enumerate}
\end{definition}

\begin{discussion}
  It takes perhaps a moment to realize that anything happened at all,
  that we defined something new here. But consider the following
  example: $X=\mb P^2$, $B$ is a rational cubic curve with a node, and
  $B'\simeq \mb P^1$ is the normalization of $B$, and $f:B'\to X$ is a
  branchdivisor whose shadow is $B$. Then the 
  pairs $(X,B)$ and $(X,B')$ can never appear as fibers in a proper family
  with connected base $S$. Indeed, $p_a(B)=1$ and $p_a(B')=0$, and the
  arithmetic genus is locally constant in proper flat families.
\end{discussion}

\begin{definition}
  A branchpair $(X,B)$ together with a morphism $f:X\to Y$ is
  \defn{stable} if $(X,\varphi_* B)$ has semi log canonical
  singularities and $K_X+\varphi_*B$ is ample over $Y$. 
\end{definition}

Finally, we define the moduli functor of stable branchpairs over a
projective scheme $Y$.

\begin{definition}
We choose a triple of positive integers 
numbers $C=$\linebreak $(C_1,C_2,C_3)$ and a positive integer $N$. We also fix a very ample
sheaf $\mc O_Y(1)$ on $Y$.  Then the basic moduli functor $M_{C,N}$ associates
to every Noetherian scheme $S$ over a base scheme the set $M_{C,N}(S)$ of
morphisms $f:X\to Y\times S$ and $\varphi_j:B_j\to X$ 
with the following properties:
\begin{enumerate}
\item $X$ and $B_j$ are flat schemes over $S$.
\item Every geometric fiber $(X,B)_s$ is a branchpair,
\item The double dual $\mc L_N(X/S)=\big( \omega_{X/S}^{\otimes N}
  \otimes \mc O_X( N\varphi_*B )\big)^{**}$ is an invertible sheaf on
  $X$, relatively ample over $Y\times S$.
\item For every geometric fiber, 
  $(L_N)_s^2=C_1$, $(L_N)_sH_s=C_2$, and $H_s^2=C_3$,  
where $\mc O_X(H)=f^*\mc O_Y(1)$. 
\end{enumerate}
\end{definition}

\begin{theorem}[Properness with normal generic fiber]
  Every family in $M_{C,N}$ over a punctured smooth curve $S\setminus
  0$ with normal $X_{\eta}$ has at most one extension, and the
  extension does exist after a ramified base change $S'\to S$.
\end{theorem}
\begin{proof}
  \emph{Existence.} The construction of the previous section gives an
  extension for the family of shadows $(X_U,\varphi_* B_U)$. The
  properness of the functor of branchvarieties \cite{AlexeevKnutson}
  applied over $X/S$ gives extensions $\varphi_j:B_j\to X$. 

  The shadow pair $(X,\varphi_* B+X_0)$ has log canonical singularities.  We
  have established that $X_0$ is $S_2$. Now by the easy direction of
  the Inversion of Adjunction 
(see, e.g. \cite[17.3]{Flips_and_Abundance})
the central fiber $(X_0,(\varphi_*B)_0)$
  has semi log canonical singularities, and so we have the required
  extension. 

  \emph{Uniqueness}. We apply Inversion of Adjunction
  \cite{Kawakita_Inversion} to the shadow pair. The conclusion is that
  $(X,\varphi_*B)$ is lc. But then it is the log canonical model of any
  resolution of singularities of any extension of
  $(X_U,\varphi_*B_U)$. Since the log canonical model is unique, the
  extension of the shadow is unique. And by the properness of the
  functor of branchvarieties \cite{AlexeevKnutson} again, the
  extensions of the branchdivisors $\varphi_j:B_j\to X$ are unique as
  well. 
\end{proof}

\begin{remark}
  One can easily see what happens when we apply our procedure in
  Example~\ref{exmp:Hassett}. The limit branchdivisor, call it $B_0'$,
  is the curve obtained from $\wB\ut_0$ by identifying two pairs of the
  points $P_i$ separately.  The morphism $B_0'\to X_0$ is 2-to-1 above
  the point $P\in X_0$ and a closed embedding away from $P$.

  Indeed, the surface $\wB$ has rational double points, and each of
  the lines in $B_0\uo\sim 2s_0$ is a $(-2)$-curve on the resolution
  of this surface. This implies that these curves are contractible.
  Let $B'$ be the projective surface obtained by contracting them.
  Thus, the central fiber $B'_0$ is obtained from $B_0\uo$ by
  identifying separately $P_1$ with $P_2$ and $P_3$ with $P_4$.

  Then the curve $(B'_0)_{\rm red}$ is nodal and has arithmetic genus
  35, the same as the generic fiber. Therefore, $B'_0 = (B'_0)_{\rm
    red}$. Hence, $B'\to C$ is a family of branchcurves. 
\end{remark}


\newcommand{\Diff}{\operatorname{Diff}}
\begin{remark}
  Examples given by J. Koll\'ar in \cite{Kollar_TwoExamples} show that
  the case of nonnormal generic fiber requires extreme care. One key
  insight from 
  \cite{Kollar_TwoExamples} is that on a properly defined non-normal
  stable pair, for every component of the double locus, the two ways
  of applying adjunction should match. 

  More precisely, let
  $(X,B)$ be a non-normal stable pair and let $C$ be
  a component of the double locus. Let $\nu:X^{\nu}\to X$ be the
  normalization, $C^{\nu}=\nu\inv(C)$, and $C^{\nu}\to C$ be the
  corresponding double cover. Then the divisor $\Diff$, computed from
  $\nu^*(K_X+B) |_{C^{\nu}} = K_{C^{\nu}}+ \Diff$, should be invariant
  under the involution.

  We also note that M.~A. van Opstall considered the case of nonnormal
  surfaces with $B=\emptyset$ in \cite{vanOpstall_Properness}.
\end{remark}





\section{Branchcycles}
\label{sec:branchcycles}

Once we have defined the branchdivisors, it is straightforward to define
branchcycles as well: the prime $k$-branchcycles of $X$ are simply
$k$-dimensional branchvarieties over $X$, and they are free generators
of an abelian group $bZ_k(X)$, resp. a free $A$-module
$bZ_k(X,A)=bZ_k(X)\otimes A$. 

\begin{definition}
  The linear (resp. algebraic) equivalence between $k$-branch\-cycles is
  generated by the following: for any family of $k$-dimensional
  branchvarieties $\phi:B\to X\times C$,
  where $C=\bP^1$ (resp. a smooth curve with two points $0,\infty$)
  the fibers $\phi_0:B_0\to X$ and $\phi_{\infty}:B_{\infty}\to X$ are
  equivalent. 

  We denote the quotient modulo the linear (resp. algebraic) equivalence 
  by $bA_k(X)$ (resp. $bB_k(X)$). 
\end{definition}

Clearly, any function which is constant in flat families of
branchvarieties descends to an invariant of $bB_k(X)$. For example,
there exists a natural homomorphism $bB_k(X)\to K_{\bullet}(X)$ to the
$K$-group of $X$ given by associating to a branchvariety $\phi:B\to X$
the class of the coherent sheaf $\phi_*\cO_B$. If we fix a very ample sheaf
$\cO_X(1)$ then we can compose this homomorphism with taking the Hilbert
polynomial to obtain a homomorphism $h:bB_k(X)\to \mb Z[t]$.

\bibliographystyle{amsalpha}

\def\cprime{$'$}

\end{document}